\documentclass[10pt,letterpaper]{article}
\usepackage{amsmath,amsfonts,amsthm,amssymb}

\theoremstyle{plain}

\newtheorem{thm}{Theorem}[section]
\newtheorem{cor}[thm]{Corollary}
\newtheorem{lem}[thm]{Lemma}

\newenvironment{exam}[1]%
{\begin{flushleft}\textbf{Example #1}.\enspace}%
{\end{flushleft}}

\newenvironment{examcont}[1]%
{\begin{flushleft}\textbf{Example #1} (continued).%
\enspace}%
{\end{flushleft}}

\newcounter{cond}

\newcommand{\real}{{\mathbb R}}
\newcommand{\bscript}{{\mathcal B}}
\newcommand{\dscript}{{\mathcal D}}
\newcommand{\escript}{{\mathcal E}}
\newcommand{\fscript}{{\mathcal F}}
\newcommand{\pscript}{{\mathcal P}}

\newcommand{\optr}{\operatorname{tr}}
\newcommand{\opsupp}{\operatorname{supp}}
\newcommand{\opker}{\operatorname{Ker}}

\newcommand{\cupdot}{\dot{\cup}}
\newcommand{\ahat}{\widehat{a}}
\newcommand{\Ahat}{\widehat{A}}
\newcommand{\bhat}{\widehat{b}}
\newcommand{\chat}{\widehat{c}}
\newcommand{\Chat}{\widehat{C}}
\newcommand{\fhat}{\widehat{f}}
\newcommand{\hhat}{\widehat{h}}
\newcommand{\rbar}{\overline{R}}

\newcommand{\ab}[1]{\left|#1\right|}
\newcommand{\brac}[1]{\left\{#1\right\}}
\newcommand{\elbows}[1]{{\left\langle#1\right\rangle}}
\newcommand{\sqbrac}[1]{\left[#1\right]}
\newcommand{\paren}[1]{\left(#1\right)}

\begin{document}

\title{UNIQUENESS AND ORDER IN\\
SEQUENTIAL EFFECT ALGEBRAS}
\author{Stan Gudder\\
Department of Mathematics\\
University of Denver\\
Denver, Colorado 80208\\
sgudder@math.du.edu\\
\and
Richard Greechie\\
Department of Mathematics\\
Louisiana Tech University\\
Ruston, Louisiana 71272\\
greechie@math.latech.edu}
\date{}
\maketitle

\begin{abstract}
A sequential effect algebra (SEA) is an effect algebra on
which a sequential product is defined. We present examples
of effect algebras that admit a unique, many and no
sequential product. Some general theorems concerning unique
sequential products are proved. We discuss sequentially
ordered SEA's in which the order is completely determined
by the sequential product. It is demonstrated that intervals
in a sequential ordered SEA admit a sequential product.
\end{abstract}

\section{Introduction} 
Quantum effects are a basic concept in foundational studies
of quantum physics \cite{blm91, bs98, da76, do66,dp00}.
Quantum effects correspond to yes-no measurements that may
be unsharp and in recent years they have been studied
within a general algebraic framework called an effect
algebra \cite{bf97,dp00,fb94,gg89,kc94}. Although effect
algebras have been useful for our understanding of quantum
theory they appear to be too general. Effect algebras only
describe one measurement connective OR (denoted by
$\oplus$) and a negation NOT (denoted by ${}'\,$). Roughly
speaking, $\oplus$ represents a parallel measurement of two
effects. However, it is important to have a mechanism for
describing series or sequential measurements of effects
(denoted by $\circ$). For this reason the authors have
introduced a structure called a sequential effect algebra
(SEA) \cite{guhist98,gg1}. A SEA is an effect algebra on
which a sequential product $\circ$ with natural properties
is defined. These properties hold in the important
case of Hilbert space effect algebras \cite{gn1,gn2}.

The present paper concentrates on uniqueness and order
properties of SEA's. We shall show that some effect algebras
admit a unique sequential product, others admit many and
still others admit none. We also present some general
results on the uniqueness of sequential products. We then
discuss a class of SEA's in which the order is completely
determined by the sequential product. We call this class
sequentially ordered SEA's. We give examples of SEA's that
are sequentially ordered and examples that are not. We
finally show that intervals in a sequentially ordered SEA
admit a sequential product. Although we review some of the
basic properties of effect algebras and SEA's we refer the
reader to our cited literature for more details.

\section{Effect Algebras} 
An \textbf{effect algebra} is an algebraic system
$(E,0,1,\oplus )$ where $0\ne1\in E$ and $\oplus$ is a
partial binary operation on $E$ satisfying:

\begin{list}
{(A\arabic{cond})}{\usecounter{cond}
\setlength{\rightmargin}{\leftmargin}}
\item If $a\oplus b$ is defined, then $b\oplus a$ is
defined and $b\oplus a=a\oplus b$.
\item If $a\oplus b$ and $(a\oplus b)\oplus c$ are defined,
then $b\oplus c$ and $a\oplus (b\oplus c)$ are defined and
$a\oplus (b\oplus c)=(a\oplus b)\oplus c$.
\item For every $a\in E$ there exists a unique $a'\in E$
such that $a\oplus a'=1$.
\item If $a\oplus 1$ is defined, then $a=0$.
\end{list}
If $a\oplus b$ is defined, we write $a\perp b$. We define
$a\le b$ if there exists a $c\in E$ such that
$a\oplus c=b$. It can be shown that $(E,\le ,{}'\,)$ is a
poset with $0\le a\le 1$ for every $a\in E$, $a''=a$, and
$a\le b$ implies $b'\le a'$ \cite{dp00,fb94}. Also,
$a\perp b$ if and only if $a\le b'$. If
$a\oplus a\oplus\cdots\oplus a\ (n\ \mathrm{summands})$ is
defined we denote this element by $na$. An element $a\in E$
is \textbf{sharp} if $a\wedge a'=0$. If every element of
$E$ is sharp, then $E$ is an \textbf{orthoalgebra}. It is
easy to show that $E$ is an orthoalgebra if and only if
$a\perp a$ implies $a=0$.

\begin{exam}{1}
For a Boolean algebra $\bscript$, define $a\perp b$ if
$a\wedge b=0$ and in this case $a\oplus b=a\vee b$. Then
$(\bscript ,0,1,\oplus )$ is an effect algebra that happens
to be an orthoalgebra. In particular, if $X\ne\emptyset$,
then $(2^X,\emptyset ,X,\oplus )$ is an effect algebra.
\end{exam}

\begin{exam}{2}
For $[0,1]\subseteq\real$, define $a\perp b$ if $a+b\le 1$
and in this case $a\oplus b=a+b$. Then
$\paren{[0,1],0,1,\oplus}$ is an effect algebra. The only
sharp elements are $0$ and $1$.
\end{exam}

\begin{exam}{3}
Let $X\ne\emptyset$ and let $\fscript\subseteq [0,1]^X$. We
call $\fscript$ a \textbf{fuzzy set system} on $X$ (i) if the
functions $0,1\in\fscript$, (ii) if $f\in\fscript$ then
$1-f\in\fscript$, (iii) if $f,g\in\fscript$ with $f+g\le 1$ then
$f+g\in\fscript$ and (iv) if $f,g\in\fscript$ then
$fg\in\fscript$. Then $(\fscript ,0,1,\oplus )$ is an
effect algebra with $f\oplus g =f+g$ whenever $f+g\le 1$.
If $\fscript =[0,1]^X$ we call $\fscript$ a \textbf{full}
fuzzy set system. The sharp elements of a fuzzy set system
$\fscript$ are the characteristic functions in $\fscript$
which can be identified with the (sharp) subsets of $X$ in
$\fscript$. Indeed, if $f\in\fscript$ is sharp, then
$f(1-f)\le f$, $1-f$ implies that $f(1-f)=0$. Hence, if
$f(x)\ne 0$, then $f(x)=1$.
\end{exam}

\begin{exam}{4}
Let $H$ be a Hilbert space and let $\escript (H)$ be the
set of self-adjoint operators on $H$ satisfying
$0\le A\le I$. For $A,B\in\escript (H)$ we define $A\perp B$
if $A+B\in\escript (H)$ and in this case $A\oplus B=A+B$.
Then $\paren{\escript (H),0,I,\oplus}$ is an effect algebra.
The elements of $\escript (H)$ are called
\textbf{quantum effects} and are important in quantum
measurement theory \cite{blm91,bs98, da76, kr83, lu83}. The
sharp elements of $\escript (H)$ are the set of projection
operators $\pscript (H)$ on $H$.
\end{exam}

\begin{exam}{5}
There are many examples of finite nonboolean effect
algebras. The simplest example is the $3$-\textbf{chain}
$C_3=\brac{0,a,1}$ where $2a=1$. Another example is the
\textbf{diamond} $D=\brac{0,a,b,1}$ where $2a=2b=1$.
\end{exam}

\begin{exam}{6}
Let $E=\omega +\omega ^*$ be the set of elements
\begin{equation*}
E=\brac{0,a,2a,\ldots ,(2a)',a',1}
\end{equation*}
By convention $0a=0$. Define $\oplus$ on $E$ by
\begin{equation*}
(ma)\oplus (na)=(m+n)a
\end{equation*}
and when $n\le m$
\begin{equation*}
(ma)'\oplus (na)=(na)\oplus (ma)'=\paren{(m-n)a}'
\end{equation*}
then $(E,0,1,\oplus )$ is an effect algebra.
\end{exam}

Let $(E_i,0_i,1_i,\oplus _i)$ be a collection of effect
algebras. One way of constructing a new effect algebra is
by taking the cartesian product $\Pi E_i$ and defining
$\Pi a_i\perp\Pi b_i$ if $a_i\perp b_i$ for every $i$ in
which case $\Pi a_i\oplus\Pi b_i=\Pi (a_i\oplus b_i)$. Then
$(\Pi E_i,\Pi 0_i,\Pi 1_i,\oplus )$ is an effect algebra.
Another way is the \textbf{horizontal sum} construction
$E=HS(E_i,i\in I)$ defined as follows. Identify all $0_i$
with a single element $0$ and all the $1_i$ with a single
element $1$. Let $E'_i=E_i\smallsetminus\brac{0_i,1_i}$,
form the disjoint union $\cupdot E_i'$ and let
$E=\brac{0,1}\cupdot E_i'$. For $a,b\in E_i$ for some
$i\in I$, if $a\perp b$ define $a\oplus b=a\oplus _ib$ and
no other orthosums are defined on $E$. Then
$(E,0,1,\oplus )$ is an effect algebra. For example
$D=HS(C_3,C_3)$. Finally, if $E$ is an effect algebra and
$b\in E$ with $b\ne 0$, then the interval
$[0,b]=\brac{a\in E\colon 0\le a\le b}$ can be organized
into an effect algebra as follows. If $c,d\in[0,b]$ and
$c\oplus d\le b$ we define $c\oplus _bd=c\oplus d$. Then
$\paren{[0,b],0,b,\oplus _b}$ is an effect algebra.

\section{Sequential Effect Algebras} 
For a binary operation $\circ$, if $a\circ b=b\circ a$ we
write $a\mid b$. A \textbf{sequential effect algebra} (SEA)
is an algebraic system $(E,0,1,\oplus ,\circ )$ where
$(E,0,1,\oplus )$ is an effect algebra and
$\circ\colon E\times E\to E$ is a binary operation
satisfying:

\begin{list}
{(S\arabic{cond})}{\usecounter{cond}
\setlength{\rightmargin}{\leftmargin}}
\item The map $b\mapsto a\circ b$ is additive for every
$a\in E$, i.e., if $b\perp c$ then $a\circ b\perp a\circ c$ and
$a\circ (b\oplus c)=a\circ b\oplus a\circ c$.
\item $1\circ a=a$ for every $a\in E$.
\item If $a\circ b=0$, then $a\mid b$.
\item If $a\mid b$, then $a\mid b'$ and
$a\circ (b\circ c)=(a\circ b)\circ c$ for every $c\in E$.
\item If $c\mid a$ and $c\mid b$, then $c\mid a\circ b$ and,
when $a\perp b$, $c\mid (a\oplus b)$.
\end{list}

We call an operation satisfying (S1)--(S5) a
\textbf{sequential product} on $E$. If $a\mid b$ for every
$a,b\in E$, then $E$ is a \textbf{commutative} SEA. Notice
that if $\circ$ is a commutative binary operation on an
effect algebra $E$, to test whether $\circ$ is a sequential
product we need only verify (S1), (S2) and
$a\circ (b\circ c)=(a\circ b)\circ c$.

Given an effect algebra $E$, does $E$ admit a sequential
product and if so is it unique? We shall show that anything
goes. There exist effect algebras that do not admit a
sequential product. There are effect algebras that admit a
unique sequential product and effect algebras that admit
many sequential products.

\begin{examcont}{1}
A Boolean algebra is a SEA under the operation $a\circ
b=a\wedge b$. It is shown in \cite{gg1} that $\circ$ is
unique.
\end{examcont}

\begin{examcont}{2}
The unit interval $[0,1]\subseteq\real$ is a SEA under the
operation $a\circ b=ab$. We shall show later that $\circ$
is unique.
\end{examcont}

\begin{examcont}{3}
A fuzzy set system $\fscript$ is a SEA under the operation
$f\circ g=fg$. We shall show later that if $\fscript$ is
full, then $\circ$ is unique.
\end{examcont}

\begin{examcont}{4}
The effect algebra $\escript (H)$ is a SEA under the
operation $A\circ B=A^{1/2}BA^{1/2}$ \cite{gn1,gn2}. This
SEA is important for quantum measurement theory
\cite{blm91,bs98} and is our first example of a
noncommutative SEA. We do not know whether $\circ$ is
unique. However, as we shall later show, $\circ$ is unique
if it satisfies some additional conditions.
\end{examcont}

\begin{examcont}{5}
The effect algebras $C_3$ and $D$ do not admit sequential
products. For $C_3$ suppose we have a sequential product
$\circ$. Then
\begin{equation*}
a=a\circ 1=a\circ (a\oplus a')=a\circ a\oplus a\circ a'
=2(a\circ a)
\end{equation*}
But there is no such element $a\circ a$ in $C_3$ which is
a contradiction. A similar demonstration holds for $D$.
More generally, it is shown in \cite{gg1} that no nonboolean
finite effect algebra admits a sequential product.
\end{examcont}

\begin{examcont}{6}
It is shown in \cite{gg1} that $E=\omega +\omega ^*$
admits a unique sequential product. For $x,y\in E$ this
sequential product is defined by
\begin{equation*}
x\circ y=
\begin{cases}0&\text{if $x=ma, y=na$}\\
[.25pc]x\wedge y
&\text{if $x=ma$, $y=(na)'$ or $x=(ma)'$, $y=na$}\\
[.25pc]\paren{(m+n)a}'&\text{if $x=(ma)'$, $y=(na)'$}
\end{cases}
\end{equation*}
\end{examcont}

Let $\dscript (H)$ be the set of density operators on $H$.
Notice that there exist $W\in\dscript (H)$ such that,
for $A\in\escript (H)$, $\optr (WA)=0$  implies $A=0$. We
call such a $W$ \textbf{faithful}. Indeed, let $x_i$ be an
orthonormal basis for $H$ and denote the one-dimensional
projection onto the span of $x_i$ by $P_i$. Then
$W=\sum\lambda _iP_i$ where $\lambda _i>0$,
$\sum\lambda _i=1$ is faithful.

\begin{exam}{7}
This is an example of an effect algebra that admits many
sequential products. Let $E_1=\escript (H)$,
$E_2=[0,1]\subseteq\real$ and $E=HS(E_1,E_2)$. Define
$\circ\colon E\times E\to E$ as follows. If $A,B\in E_1$
then $A\circ B=A^{1/2}BA^{1/2}$; if $a,b\in E_2$ then
$a\circ b=ab$; if $A\in E_1$, $a\in E_2$ then
$A\circ a=aA$; and if $a\in E_2$, $A\in E_1$ then
$a\circ A=a\optr (WA)$ where $W\in\dscript (H)$ is fixed and
faithful. It is shown in \cite{gg1} that $\circ$ is
sequential product on $E$. Notice that different faithful
$W\in\dscript (H)$ give different sequential products.
\end{exam}

Let $a$ be a sharp element of an effect algebra $E$.
Suppose we view $a$ in a larger context by enlarging $E$ to
an effect algebra $F$. Since $a$ may not be sharp as a member
of $F$ we say that sharpness is \textbf{contextual} in effect
algebras. Physically we would not expect sharpness to be
contextual and this is an unfortunate property for effect
algebras. The next result shows that this unfortunate
property holds for any effect algebra that contains a
nontrivial sharp element. The result also shows that
sharpness is noncontextual for SEA's. We denote the set of
sharp elements of an effect algebra $E$ by
$E_S$. As usual, an \textbf{embedding} for an effect algebra
is a monomorphism \cite{dp00,fb94}.

\begin{thm}      
\label{thm31}
{\rm (i)}\enspace If $E$ is an effect algebra, then there
exists an effect algebra $F$ and an effect algebra embedding
$\phi\colon E \to F$ such that $F_S=\brac{ 0,1}$.
{\rm (ii)}\enspace If $E$ and $F$ are SEA's and
$\phi\colon E\to F$ is a SEA embedding, then $\phi (a)\in
F_S$ if and only if $a\in E_S$.
\end{thm}

\begin{proof}
(i)\enspace Let $G$ be a nontrivial abelian partially
ordered group. Define
\begin{equation*}
E_G=\paren{E\smallsetminus\brac{0,1}}
\times G\cup\brac{(0,g),(1,-g)\colon g\ge 0}
\end{equation*}
and for $(a,g),(b,h)\in E_G$ define
$(a,g)\oplus (b,h)=(a\oplus b,g+h)$ provided that $a\perp b$
and $(a\oplus b,g+h)\in E_G$. Letting $\mathbf{0}=(0,0)$,
$\mathbf{1}=(1,0)$, we shall show that
$(E_B,\mathbf{0},\mathbf{1},\oplus )$ is an effect
algebra. It is clear that (A1) (commutativity) holds.
Defining $(a,g)'=(a',-g)$ it is easy to check that (A3)
holds. To verify (A4), suppose that $(a,g)\oplus\mathbf{1}$
is defined. It follows that $a=0$ and $(1,g)\in E_G$. Since
$(0,g)\in E_G$, it follows that $g=0$. Hence,
$(a,g)=\mathbf{0}$.

To verify (A2) (associativity), assume that
$(a,g),(b,h),(c,k)\in E_G$ with $(b,h)\oplus (c,k)\in E_G$
and
\begin{equation*}
(a,g)\oplus\sqbrac{(b,h)\oplus (c,k)}\in E_G\quad (*)
\end{equation*}
Then
\begin{align*}
(a,g)\oplus\sqbrac{(b,h)\oplus (c,k)}
&=\paren{a\oplus (b\oplus c),g+(h+k)}\\
&=\paren{(a\oplus b)\oplus c,(g+h)+k}\\
&=(a\oplus b, g+h)\oplus (c,k)
\end{align*}
provided that $(a\oplus b,g+h)\in E_G$. Noting that this
holds for all $g,h\in G$ when $a\oplus b\notin\brac{0,1}$,
we need only consider the cases $a\oplus b\in\brac{0,1}$. If
$a\oplus b=0$, then $a=b=0$. Hence, $g,h\ge 0$ so that
$(a\oplus b,g+h)\in E_G$. If $a\oplus b=1$, then $c=0$ so
that $k\ge 0$. Also, by (*), we have $(g+h)+k\le 0$ so that
$g+h\le -k\le 0$. Hence, $(a\oplus b,g+h)\in E_G$ so
$(E_G,\mathbf{0},\mathbf{1},\oplus )$ is an effect algebra.

Note that $(a,g)\le (b,h)$ in $E_G$ if and only if $a\le b$
and $(b\ominus a,h-g)\in E_G$. We thus have either $a< b$ or
$a=b$ and $g\le h$ which is the lexicographic order on
$E_G$. Define $\phi\colon E\to E_G$ by $\phi (a)=(a,0)$.
Clearly $\phi$ is an effect algebra embedding of $E$ into
$E_G$. To show that $(E_G)_S=\brac{\mathbf{0},\mathbf{1}}$
suppose $(a,g)\in E_G$ with $a\ne 0,1$. Then for $h\in G$
with $h>0$ we have $(0,h)\le (a,g)$, $a',-g)$. Hence,
$(a,g)\notin (E_G)_S$. For the case $(0,g)\in E_G$ with
$g>0$ we have $(0,g)<(1,-g)$ and for the case $(1,g)\in E_G$
with $g<0$ we have $(0,-g)<(1,g)$.

(ii)\enspace If $a\in E_S$, then
$\phi (a)\circ\phi (a)=\phi (a\circ a)=\phi (a)$ so that
$\phi (a)\in F_S$. Conversely, if $a\notin E_S$ then there
exists a $b\in E$ such that $0<b\le a,a'$. But then
$0<\phi (b)\le\phi (a),\phi (a)'$ so that $\phi (a)\notin
F_S$.
\end{proof}

The second part of the proof of Theorem~\ref{thm31}(ii)
shows that fuzziness is noncontextual in an effect algebra
$E$. That is, if $\phi\colon E\to F$ is an effect algebra
embedding and $a\notin E_S$ then $\phi (a)\notin F_S$. The
construction of the effect algebra $E_G$ in the proof of
Theorem~\ref{thm31}(i) is of interest in its own right. The
elements of the form $(0,g)\in E_G$ act like infinitesimals.
In the special case where $E=\brac{0,1}$ is trivial and $G$
is the integers $Z$, we have that $E_Z$ is isomorphic to
$\omega +\omega ^*$.

\section{Results} 
\begin{thm} 
\label{thm41}
There is a unique sequential product on the effect algebra
$[0,1]\subseteq\real$.
\end{thm}
\begin{proof}
Let $\circ$ be a sequential product on $[0,1]$. Then for
any integer $n\ge 1$ and $a\in [0,1]$ we have
\begin{equation*}
a=a\circ 1
=a\circ\paren{\frac{1}{n}\oplus\cdots\oplus\frac{1}{n}}
=n\paren{a\circ\frac{1}{n}}
\end{equation*}
so that $a\circ\frac{1}{n}=\frac{1}{n}a$. Also, for any
integer $1\le m\le n$ we have
\begin{equation*}
a\circ\frac{m}{n}
=a\circ\paren{\frac{1}{n}\oplus\cdots\oplus\frac{1}{n}}
=m\paren{a\circ\frac{1}{n}}=\frac{m}{n}\,a
\end{equation*}
Hence, for any rational number $r\in Q\cap [0,1]$ we have
$a\circ r=ar$. Now let $b\in [0,1]$ be irrational. If
$r\in Q\cap [0,1]$ and $b<r$ then by additivity we obtain
\begin{equation*}
a\circ b\le a\circ r=ar
\end{equation*}
Similarly, if $r\in Q\cap [0,1]$ and $b>r$, then
\begin{equation*}
a\circ b\ge a\circ r=ar
\end{equation*}
Since $Q\cap [0,1]$ is dense in $[0,1]$, we obtain
$a\circ b=ab$.
\end{proof}

The next theorem formalizes the obvious observation that if
two effect algebras are isomorphic and one admits an
operation satisfying some special conditions, then so does
the other.

\begin{thm} 
\label{thm42}
Let $E,F$ be effect algebras and let $\phi\colon E\to F$ be
an effect algebra isomorphism {\rm\cite{dp00,fb94}}. If
$\circ$ is a sequential product on $E$, then
$a*b=\phi\sqbrac{\phi ^{-1}(a)\circ\phi ^{-1}(b)}$ is a
sequential product on $F$. Moreover, $(E,\circ )$ and
$(F,*)$ are SEA isomorphic {\rm\cite{gg1}}.
\end{thm}
\begin{proof}
The proof is a straightforward verification.
\end{proof}

\begin{cor} 
\label{cor43}
If $E$ and $F$ are isomorphic effect algebras and $E$
admits a unique sequential product $\circ$, then $F$ admits
a unique sequential product.
\end{cor}
\begin{proof}
By Theorem~\ref{thm42}, $F$ admits a sequential product
$*$. Let $\phi\colon E\to F$ be an effect algebra
isomorphism and define $\bullet\colon E\times E\to E$ by
$a\bullet b=\phi ^{-1}\sqbrac{\phi (a)*\phi (b)}$. By
Theorem~\ref{thm42}, $\bullet$ is a sequential product on
$E$ so that $a\bullet b=a\circ b$. Hence,
\begin{equation*}
\phi (a\circ b)=\phi (a\bullet b)=\phi (a)*\phi (b)
\end{equation*}
thus, for every $c,d\in F$ we have
\begin{equation*}
c*d=\phi \sqbrac{\phi ^{-1}(c)}*\phi\sqbrac{\phi ^{-1}(d)}
=\phi\sqbrac{\phi ^{-1}(c)\circ\phi ^{-1}(d)}
\end{equation*}
It follows that $*$ is unique.
\end{proof}

\begin{cor} 
\label{cor44}
If an effect algebra $E$ admits a unique sequential product
$\circ$, then any effect algebra automorphism
$\phi\colon E\to E$ is a SEA automorphism.
\end{cor}
\begin{proof}
By Theorem~\ref{thm42},
$a*b=\phi ^{-1}\sqbrac{\phi (a)\circ\phi (b)}$ is a
sequential product on $E$. Hence, $a*b=a\circ b$ and the
result follows.
\end{proof}

We say that $a,b\in E$ \textbf{coexist} if there exist
$c,d,e\in E$ such that $c\oplus d\oplus e$ is defined and
$a=c\oplus d$, $b=c\oplus e$.

\begin{lem} 
\label{lem45}
Let $\circ$ and $*$ be sequential products on an effect
algebra $E$ and let $b\in E_S$. If $a\circ b=b\circ a$,
then $a*b=b*a$.
\end{lem}
\begin{proof}
It is shown in \cite{gg1} that $a\circ b=b\circ a$ if and
only if $a$ and $b$ coexist. But coexistence is independent
of the sequential product.
\end{proof}

\begin{thm} 
\label{thm46}
Let $(E_i,0_i,1_i,\oplus _i,\circ _i)$ be SEA's, $i\in I$.
Then $E=\Pi E_i$ admits a unique sequential product if and
only if each $E_i$, $i\in I$, admits a unique sequential
product.
\end{thm}
\begin{proof}
If $E_j$ admits two sequential products for some $j\in I$,
then clearly $\Pi E_i$ admits at least two sequential
products. Conversely, suppose $E_i$, $i\in I$, admits a
unique sequential product $\circ _i$ and let $*$ be a
sequential product on $\Pi E_i$. For $j\in I$, let
$f_j\in\Pi E_i$ be define by
\begin{equation*}
f_j(i)=
\begin{cases}1_j&\text{if $i=j$}\\
[.25pc]0_i&\text{if $i\ne j$}
\end{cases}
\end{equation*}
Clearly, $f_j\in E_S$. Let $\circ$ be the sequential
product on $E$ given by $(f\circ g)(i)=f(i)\circ _ig(i)$,
$i\in I$, for any $f,g\in E$. Since
$f\circ f_j=f_j\circ f$, by Lemma~\ref{lem45},
$f*f_j=f_j*f$ for any $f\in E$. It follows from
Theorem~3.4 \cite{gg1} that $f*f_j=f\wedge f_j$. Hence,
\begin{equation*}
f*f_j(i)=
\begin{cases}0_i&\text{if $i\ne j$}\\
[.25pc]f(j)&\text{if $i=j$}
\end{cases}
\end{equation*}
For any $f,g\in E$ we have
\begin{align*}
(f*g)*f_j&=f_j*(f*g)=(f_j*f)*g
=\sqbrac{(f*f_j)*f_j}*g\\
&=(f*f_j)*(f_j*g)
\end{align*}
Now $[0,f_j]\subseteq E$ is an effect algebra with greatest
element $f_j$ and $\phi\colon E_j\to[0,f_j]$ given by
\begin{equation*}
\sqbrac{\phi (a)}(i)=
\begin{cases}0_i&\text{if $i\ne j$}\\
[.25pc]a&\text{if $i=j$}
\end{cases}
\end{equation*}
is an effect algebra isomorphism. Since $E_j$ admits a
unique sequential product, by Corollary~\ref{cor43},
$[0,f_j]$ admits a unique sequential product. Hence,
\begin{align*}
(f*g)(j)&=\sqbrac{(f*g)*f_j}(j)
=\sqbrac{(f*f_j)*(g*f_j)}(j)\\
&=\sqbrac{(f*f_j)\circ (g*f_j)}(j)
=(f*f_j)(j)\circ _j(g*f_j)(j)\\
&=f(j)\circ _jg(j)=(f\circ g)(j)
\end{align*}
Hence, $\circ$ is the unique sequential product on $E$.
\end{proof}

\begin{cor} 
\label{cor47}
A full fuzzy set system $E=[0,1]^X$ admits a unique
sequential product.
\end{cor}
\begin{proof}
This follows from Theorem~\ref{thm46} because, by
Theorem~\ref{thm41}, $[0,1]$ admits a unique sequential
product
\end{proof}

The next result characterizes the only known sequential product on $\escript (H)$.
For $x\in H$ with $\|x\|=1$, $P_x$
denoted the one-dimensional projection onto the span of $x$.

\begin{thm} 
\label{thm48}
Let
$\circ\colon\escript (H)\times\escript (H)\to\escript (H)$
be a binary operation. Then $A\circ B=A^{1/2}BA^{1/2}$ for
every $A,B\in\escript (H)$ if and only if the following
conditions are satisfied:
{\rm (1)}\enspace $B\mapsto A\circ B$ is $\sigma$-additive in
the strong operator topology for every $A\in\escript (H)$;
{\rm (2)}\enspace $(\lambda A)\circ B=\lambda (A\circ B)$
for every $\lambda\in [0,1]$;
{\rm (3)}\enspace there exists a Borel function
$f\colon [0,1]\to[0,1]$ such that $f(1)=1$ and
$\elbows{A\circ P_xy,y}=\ab{\elbows{f(A)x,y}}^2$ for every
$A\in\escript (H)$, $x,y\in H$ with $\|x\|=\|y\|=1$.
\end{thm}
\begin{proof}
We have already observed that $A\circ B=A^{1/2}BA^{1/2}$ is
a sequential operation on $\escript (H)$; that these
properties hold for this operation is straightforward, with
$f(\lambda )=\lambda ^{1/2}$ in (3). To prove the converse,
assume the conditions and observe that, by (2) and (3), for
every $\lambda\in [0,1]$ we have
\begin{align*}
\ab{\elbows{f(\lambda A)x,y}}^2
&=\elbows{(\lambda A)\circ P_xy,y}
=\lambda\elbows{A\circ P_xy,y}
=\lambda\ab{\elbows{f(A)x,y}}^2\\
&=\ab{\elbows{\lambda ^{1/2}f(A)x,y}}^2
\end{align*}
Letting $y=x$ gives
$\elbows{f(\lambda A)x,x}=\elbows{\lambda ^{1/2}f(A)x,x}$
for every $x\in H$ with $\|x\|=1$. Hence,
$f(\lambda A)=\lambda ^{1/2}f(A)$. Letting $A=I$ gives
\begin{equation*}
f(\lambda )I=f(\lambda I)=\lambda ^{1/2}f(I)
=\lambda ^{1/2}I
\end{equation*}
so that $f(\lambda )=\lambda ^{1/2}$. Thus,
\begin{equation*}
\elbows{A\circ P_xy,y}=\ab{\elbows{A^{1/2}x,y}}^2
=\elbows{A^{1/2}P_xA^{1/2}y,y}
\end{equation*}
for every $y\in H$ with $\|y\|=1$. It follows that
$A\circ P_x=A^{1/2}P_xA^{1/2}$. By (1) we have
$A\circ P=A^{1/2}PA^{1/2}$ for every $P\in\pscript (H)$. As
in the proof of Theorem~\ref{thm41} we have
$A\circ (\lambda B)=\lambda A\circ B$ for every
$\lambda\in [0,1]$. Hence, by (1) we conclude that
$A\circ B=A^{1/2}BA^{1/2}$ for every $B\in\escript (H)$ with
finite spectrum. Since any $B\in\escript (H)$ is the strong
limit of an increasing sequence of $B_i\in\escript (H)$
with finite spectra, it follows from (1) that
$A\circ B=A^{1/2}BA^{1/2}$ for every $B\in\escript (H)$.
\end{proof}

Since the sequential product $A\circ B=A^{1/2}BA^{1/2}$ is
the only known sequential product in $\escript (H)$, we
shall refer to it as the standard sequential product on
$\escript (H)$. When we refer to  $\escript (H)$ as a SEA,
without reference to a specific sequential product, we mean
with respect to the standard sequential product.

\section{Sequentially Ordered SEA's} 
An effect algebra $E$ is \textbf{sharply dominating} if for
every $a\in E$ there exists a least element $\ahat\in E_S$
such that $a\le\ahat$ \cite{gusharp98}. A sharply
dominating SEA is \textbf{sequentially ordered} if
(1)\enspace $a\le b$ implies that there exists a $c\in E$
such that $a=b\circ c$ and
(2)\enspace if $c\circ a\le c\circ b$ then
$\chat\circ a\le\chat\circ b$. Notice that the converses of
(1) and (2) hold for any SEA. Condition (1) states that the
order is completely determined by $\circ$ ($a\le b$ if and
only if $a=b\circ c$ for some $c$). This is similar to
order being completely determined by $\oplus$ ($a\le b$ if
and only if $a\oplus c=b$ for some $c$). It is easy to
check that Boolean algebras and $[0,1]\subseteq\real$ are
sequentially ordered.

\begin{examcont}{3}
Let $E=[0,1]^X$ be a full fuzzy set system. For $f\in E$ let
\begin{equation*}
\opsupp (f)=\brac{x\in X\colon f(x)\ne 0}
\end{equation*}
and define $\fhat$ to be the characteristic function on
$\opsupp (f)$. Then $\fhat$ is the least sharp element that
dominates $f$ so $E$ is sharply dominating. If $f\le g$
then $f=gh$ where
\begin{equation*}
h(x)=
\begin{cases}f(x)/g(x)&\text{if $g(x)\ne 0$}\\
[.25pc]0
&\text{if $g(x)=0$}
\end{cases}
\end{equation*}
Hence, $E$ satisfies Condition~(1). If $hf\le hg$, then
$f(x)\le g(x)$ for all $x\in\opsupp (h)$ so
$\hhat f\le\hhat g$. Hence, $E$ satisfies Condition~(2) so
$E$ is sequentially ordered.

Now let $F$ be the fuzzy set system of all polynomial
functions $f\colon [0,1]\to[0,1]$. Then $F_S=\brac{0,1}$
and $F$ is sharply dominating. But $F$ does not satisfy (1)
so $F$ is not sequentially ordered. Indeed, the functions
$f(x)=\frac{1}{2} x$ and $g(x)=\frac{1}{2}+\frac{1}{2}x$
are in $F$ and $f\le g$. Suppose there exists an $h\in F$
such that $f=gh$. Then $h(x)=x/(x+1)$ on $[0,1]$ but
$h\notin F$ which is a contradiction. However, $F$ does
satisfy (2). Indeed, if $hf\le hg$ then $f(x)\le g(x)$ for
all $x\in\opsupp (h)$. But if $h\ne 0$, then $h(x)=0$ for
only a finite number of points $x_i\in [0,1]$,
$i=1,\ldots ,n$. If $f(x_i)>g(x_i)$ for some $i$, then by
continuity $f(x)>g(x)$ in a neighborhood of $x_i$ which is
a contradiction. Hence,
\begin{equation*}
\hhat f=f\le g=\hhat g
\end{equation*}
\end{examcont}

\begin{examcont}{4}
It is well known that $\escript (H)$ is sharply dominating
\cite{gusharp98}. We shall show in the next theorem that
$\escript (H)$ is sequentially ordered.
\end{examcont}

\begin{examcont}{6}
The SEA $\omega +\omega ^*$ is sharply dominating with
$(\omega +\omega ^*)_S=\brac{0,1}$. However,
$\omega +\omega ^*$ is not sequentially ordered. Indeed,
$a\le 2a$ but there is no $c\in\omega +\omega ^*$ such that
$a=(2a)\circ c$. Also, (2) does not hold because
$a\circ 2a=0=a\circ a$ but
\begin{equation*}
\ahat\circ 2a=2a\not\le a=\ahat\circ a
\end{equation*}
\end{examcont}

\begin{examcont}{7}
It is easy to check that $E=HS\paren{\escript (H),[0,1]}$
is sharply dominating and satisfies (1). However, $E$ does
not satisfy (2). Indeed, if $a\in (0,1)$ then
$a\circ A\le a\circ B$ if and only if
$\optr (WA)\le\optr (WB)$ but this does not imply that
\begin{equation*}
\ahat\circ A=A\le B=\ahat\circ B
\end{equation*}
This observation together with the second part of Example~3
shows that Conditions~(1) and (2) are logically independent.
\end{examcont}

\begin{thm} 
\label{thm51}
The SEA $\escript (H)$ is sequentially ordered.
\end{thm}
\begin{proof}
For $A\in\escript (H)$ let $P_A$ be the projection onto the
closure of the range $\rbar (A)$ of $A$. Then
$P_S\in\pscript (H)=\escript (H)_S$ and it is easy to see
that $P_A$ is the least sharp element satisfying
$A\le P_A$. Hence, $P_A=\Ahat$ and $\escript (H)$ is sharply
dominating. We now show that if $C\in\escript (H)$, then
$\rbar (C)=\rbar (C^{1/2})$. If $Cx=0$, then
\begin{equation*}
\elbows{C^{1/2}x,C^{1/2}x}=\elbows{Cx,x}=0
\end{equation*}
so that $C^{1/2}x=0$. Hence,
$\opker (C)\subseteq\opker (C^{1/2})$. Conversely, if
$C^{1/2}x=0$ then $Cx=0$ so
$\opker (C^{1/2})\subseteq\opker (C)$. Hence,
$\opker (C)=\opker (C^{1/2})$ and we have
\begin{equation*}
\rbar (C)=\opker (C)^{\perp}=\opker (C^{1/2})^{\perp}
=\rbar (C^{1/2})
\end{equation*}
Now suppose that $C\circ A\le C\circ B$. Then for any
$x\in H$ we have
\begin{equation*}
\elbows{AC^{1/2}x,C^{1/2}x}
\subseteq\elbows{BC^{1/2}x,C^{1/2}x}
\end{equation*}
Hence,$\elbows{Ay,y}\le\elbows{By,y}$ for any
$y\in\rbar (C^{1/2})=\rbar (C)$. We conclude that
\begin{equation*}
\elbows{AP_Cx,P_Cx}\le\elbows{BP_Cx,P_Cx}
\end{equation*}
for any $x\in H$. Hence,
\begin{equation*}
\Chat\circ A=P_CAP_C\le P_CBP_C=\Chat\circ B
\end{equation*}
so Condition~(2) holds. To verify Condition~(1) suppose
that $A\le B$. Then $A^{1/2}A^{1/2}\le B^{1/2}B^{1/2}$ and
it follows from \cite{do66} that there exists a bounded
linear operator $T$ on $H$ such that $\|T\|\le 1$ and
$A^{1/2}=B^{1/2}T$. Letting $C=TT^*$, we see that $C\ge 0$.
Moreover, for any $x\in H$ we have
\begin{equation*}
\elbows{Cx,x}=\elbows{T^*x,T^*x}=\|T^*x\|^2
\le\|T^*\|^2\|x\|^2\le\|x\|^2=\elbows{x,x}
\end{equation*}
so that $C\in\escript (H)$. Hence,\medskip

\hskip 3pc
$A=A^{1/2}(A^{1/2})^*=B^{1/2}TT^*B^{1/2}=B^{1/2}CB^{1/2}
=B\circ C$
\end{proof}

\begin{thm} 
\label{thm52}
Let $E$ be a sequentially ordered SEA. For $a,b\in E$ with
$a\le b$ there exists a unique $c\in E$ such that
$c\le\bhat$ and $a=b\circ c$.
\end{thm}
\begin{proof}
By Condition~(1) there is a $d\in E$ such that
$a=b\circ d$. Letting $c=\bhat\circ d$ we have $c\le\bhat$
and
\begin{equation*}
a=b\circ d=(b\circ\bhat )\circ d=b\circ (\bhat\circ d)
=b\circ c
\end{equation*}
For uniqueness, suppose that $c_1\le\bhat$ and $a=b\circ
c_1$. Then $b\circ c_1=b\circ c$ and applying Condition~(2)
we have\medskip

\hskip 10pc $c_1=\bhat\circ c_1=\bhat\circ c=c$
\end{proof}
We denote the unique element $c$ in Theorem~\ref{thm52} by
$c=a/b$ and call $c$ the \textbf{sequential quotient} of $a$
over $b$. Thus, ${}/$ is a partial binary operation on $E$
with domain $\brac{(a,b)\colon a\le b}$.

\begin{cor} 
\label{cor53}
Let $E$ be a sequentially ordered SEA.
{\rm (i)}\enspace For every $a,b,c\in E$ there exists a
unique $d\in E$ such that $d\le (a\circ b)^\wedge$ and
$a\circ (b\circ c)=(a\circ b)\circ d$.
{\rm (ii)}\enspace $a\le b$ if and only if there exists a
unique $d\in E$ such that $d\ge (\bhat )'$ and
$a\oplus b\circ d=b$.
\end{cor}
\begin{proof}
(i)\enspace Since $b\circ c\le b$ we have
$a\circ (b\circ c)\le a\circ b$. By Theorem~\ref{thm52}
there exists a unique $d\in E$ such that
$d\le (a\circ b)^\wedge$ and
$a\circ (b\circ c)=(a\circ b)\circ d$.
(ii)\enspace If $a\le b$, then by Theorem~\ref{thm52} there
exists a unique $c\in E$ (namely, $c=a/b$) such that
$c\le\bhat$ and $a=b\circ c$. Hence, $c'\ge (\bhat )'$ and
\begin{equation*}
b=b\circ c\oplus b\circ c'=a\oplus b\circ c'
\end{equation*}
with the uniqueness of $c'$ following from the uniqueness of
$c$. The converse is clear.
\end{proof}

The proof of the next lemma is straightforward.

\begin{lem} 
Let $E$ be a sequentially ordered SEA and let $a\in E$.
{\rm (i)}\enspace $a/a=\ahat$.
{\rm (ii)}\enspace  $a\in E_S$ if and only if $a/a=a$.
{\rm (iii)}\enspace Let $b\in E_S$. If $a\le b$
then $a/b=a$ and if $b\le a$ then $b/a=b$. In particular,
$a/1=a$ and $0/a=0$ for every $a\in E$.
{\rm (iv)}\enspace If $n\ge 1$, then $a^{n+m}/a^m=a^n$.
\end{lem}

If $a\le b$ then the unique $c$ such that $a\oplus c=b$ is
denoted by $b\ominus a$.

\begin{thm} 
\label{thm55}
Let $E$ be a sequentially ordered SEA with $a,b,c \in E$.
{\rm (i)}\enspace If $a\le b$ then
$(b\ominus a)/b=\bhat\circ (a/b)'$.
{\rm (ii)}\enspace $(a\circ b)/a=\ahat\circ b$.
{\rm (iii)}\enspace If $a\le b\le c$ then $a/c\le b/c$.
{\rm (iv)}\enspace Let $a,b\le c$. Then  $a\oplus b\le c$ iff
$(a/c)\perp (b/c)$, and in this case
$(a\oplus b)/c=a/c\oplus b/c$.
{\rm (v)}\enspace
$a/(a\oplus b)=(a\oplus b)^\wedge
\circ\sqbrac{b/(a\oplus b)}'$.
{\rm (vi)}\enspace If $a\le b$ and $b\mid (a/b)$ then
$b\mid a$.
\end{thm}
\begin{proof}
(i)\enspace For $a\le b$ we have
\begin{equation*}
b\circ\sqbrac{\bhat\circ (a/b)'}=b\circ (a/b)'
=b\circ (1\ominus a/b)=b\ominus b\circ (a/b)=b\ominus a
\end{equation*}
and $\bhat\circ (a/b)'\le \bhat$.
(ii)\enspace Since $a\circ b=a\circ (\ahat\circ b)$ and
$\ahat\circ b\le\ahat$ we have $(a\circ b)/a=\ahat\circ b$.
(iii)\enspace Since $a=c\circ (a/c)$ and $b=c\circ (b/c)$
we have $c\circ (a/c)\le c\circ (b/c)$. Applying
Condition~(2) gives
\begin{equation*}
a/c=\chat\circ (a/c)\le\chat\circ (b/c)=b/c
\end{equation*}
(iv)\enspace If $a\oplus b\le c$, then $a,b\le c$ and both
$a/c$ and $b/c$ are defined. Since $a\le c\ominus b$, by (i)
and (iii) we have
\begin{equation*}
a/c\le (c\ominus b)/c=\chat\circ (b/c)'=\chat\ominus (b/c)
\end{equation*}
Hence, $a/c\oplus b/c$ is defined. Now
\begin{equation*}
a\oplus b=c\circ (a/c)\oplus c\circ (b/c)
=c\circ (a/c\oplus b/c)
\end{equation*}
Since, by Lemma 4.2 of [13], $a/c\oplus b/c\le\chat$ we have
$(a\oplus b)/c=a/c\oplus b/c$. If $(a/c)\perp (b/c)$ then,
since
\begin{equation*}
c\circ (a/c\oplus b/c)=c\circ (a/c)\oplus c\circ (b/c)
=a\oplus b
\end{equation*}
we have $a\oplus b\le c$. The result now follows from the
above.

\noindent (v)\enspace This follows from (i).

\noindent (vi)\enspace This follows because
$b\circ a=b\circ\sqbrac{b\circ (a/b)}
=\sqbrac{b\circ (a/b)}\circ b=a\circ b$
\end{proof}

The next result shows that the condition in
Theorem~\ref{thm55}(ii) characterizes the sequential
quotient.

\begin{lem} 
\label{lem56}
Let $E$ be a sequentially ordered SEA. If $/\!/$ is a
partial binary operation on $E$ with domain
$\brac{(a,b)\colon a\le b}$ and
$(a\circ b)/\!/a=\ahat\circ b$ for every $a,b\in E$, then
$/\!/$ and ${}/$ coincide.
\end{lem}
\begin{proof}
If $a\le b$ then $a=b\circ c$ where $c\le\bhat$. Hence
\smallskip

\hskip 8pc $a/\!/b=b\circ c/\!/b=\bhat\circ c=c=a/b$
\end{proof}

Let $E$ be an effect algebra and let $b\in E$ with
$b\ne 0$. Then we have seen in Section~2 that
$\paren{[0,b],0,b,\oplus _b}$ is an effect algebra. If $E$
is also a SEA, does $[0,b]$ admit a sequential product? If
$b\in E_S$ the answer is yes. Just restrict $\circ$ to
$[0,b]$. In this case, $b\circ a=a$ for all $a\in [0,b]$
and the other axioms are easily verified so that
$\paren{[0,b],0,b,\oplus _b,\circ}$ is a SEA. In general,
the answer is no. For example, in $\omega +\omega ^*$ the
interval $[0,2a]=\brac{0,a,2a}$ is isomorphic to $C_3$ so
$[0,2a]$ does not admit a sequential product. We now show
that the answer is positive if $E$ is sequentially ordered.

Let $E$ be a sequentially ordered SEA and let $b\in E$ with
$b\ne 0$. Define $\phi _b\colon [0,b]\to [0,\bhat ]$ by
$\phi _b(a)=a/b$.

\begin{lem} 
\label{lem57}
The map $\phi _b\colon [0,b]\to [0,\bhat ]$ is an effect
algebra isomorphism.
\end{lem}
\begin{proof}
By Theorem~\ref{thm55}(iv), if $a,c\in [0,b]$ and
$a\oplus c\le b$ then
\begin{equation*}
\phi _b (a\oplus c)=(a\oplus c)/b=a/b\oplus c/b
=\phi _b(a)\oplus \phi _b(c)
\end{equation*}
Hence, $\phi _b$ is additive. Also, $\phi _b(b)=\bhat$ so
$\phi _b$ is a morphism \cite{dp00,fb94}. If $\phi
_b(a)\perp\phi _b(c)$ then by Theorem~\ref{thm55}(v),
$a\oplus c\le b$ so $a\perp c$. Thus, $\phi _b$ is a
monomorphism \cite{dp00,fb94}. If $c\in [0,\bhat ]$,
letting $a=b\circ c$ we have $a\in [0,b]$ and
$\phi _b(a)=a/b=c$. Hence, $\phi _b$ is surjective so
$\phi _b$ is an effect algebra isomorphism.
\end{proof}

\begin{thm} 
\label{thm58}
Let $E$ be a sequentially ordered SEA and let $b\in E$ with
$b\ne 0$.
{\rm (i)}\enspace There exists a unique sequential product
$\circ _b$ on $[0,b]$ such that
\begin{equation*}
(a\circ _bc)/b=(a/b)\circ (c/b)
\end{equation*}
{\rm (ii)}\enspace Employing this sequential product on
$[0,b]$,
$\phi _b\colon [0,b]\to [0,\bhat ]$ becomes a SEA
isomorphism.
\end{thm}
\begin{proof}
(i)\enspace Uniqueness follows from
$a\circ _bc=b\circ\sqbrac{(a/b)\circ (c/b)}$. By
Lemma~\ref{lem57}, $\phi _b\colon [0,b]\to [0,\bhat ]$ is an
effect algebra isomorphism. Letting $\psi =\phi _b^{-1}$ we
conclude, using  Theorem~\ref{thm55} part (ii),
that $\psi\colon [0,\bhat ]\to [0,b]$ is an effect
algebra isomorphism given by $\psi (a)=b\circ a$. Applying
Theorem~\ref{thm42}, for $a,c\in [0,b]$ we have
\begin{equation*}
a*c=\psi\sqbrac{\psi ^{-1}(a)\circ\psi ^{-1}(b)}
=b\circ\sqbrac{\phi _b(a)\circ\phi _b(c)}
=b\circ\sqbrac{(a/b)\circ (c/b)}
\end{equation*}
is a sequential product on $[a,b]$.
(ii)\enspace This follows from Theorem~\ref{thm42}.
\end{proof}

\end{document}